%% file: adatyp2h.tex
\documentclass[12pt]{article}
\usepackage[reqno]{amsmath}
\usepackage{amssymb, theorem, enumerate}
\usepackage{graphicx}

\input Lmacs

\expandafter\ifx\csname pdfoptionalwaysusepdfpagebox\endcsname\relax\else
\fi

\theoremstyle{change}
{\theorembodyfont{\slshape}
\newtheorem{theorem}{Theorem.}[section]
\newtheorem{lemma}[theorem]{Lemma.}
\newtheorem{corollary}[theorem]{Corollary.}}
\theorembodyfont{\rmfamily}

\newcommand\lref[1]{Lemma~\ref{lem:#1}}
\newcommand\tref[1]{Theorem~\ref{thm:#1}}
\newcommand\cref[1]{Corollary~\ref{cor:#1}}

\def\proof{\noindent{\textsl{Proof. }}}
\def\sqr#1#2{{\vbox{\hrule height.#2pt
    \hbox{\vrule width.#2pt height#1pt \kern#1pt
        \vrule width.#2pt}\hrule height.#2pt}}}
\def\eqed{\sqr53}
\def\qed{%
    \ifmmode\eqno\eqed
    \else\nobreak\ \hfill\eqed\medbreak\fi}

\newcommand\yy[2]{Y_{{#1},{#2}}} 
\newcommand\yab{\yy ab}
\newcommand\yij{\yy ij}
\newcommand\ypab{Y'_{a,b}}
\newcommand\dist{\opk{dist}}

\def\ix{x^{-1}}
\def\iy{y^{-1}}
\def\iz{z^{-1}}

\def\nom#1{\cN_{#1}}
\def\snv{^{(-)}}
\def\snt{^{(-)T}}

\def\thm{\Th(M)}
\def\thmw{\Th_W(M)}
\def\nomw{\nom{W}}
\def\snt{^{(-)T}}
\def\snv{^{(-)}}

\title{Type-II Matrices and Combinatorial Structures}
\author{Ada Chan and Chris Godsil}
\begin{document}
\maketitle

\begin{abstract}
Type-II matrices are a class of matrices used by Jones in his work on spin models.  
In this paper we show that type-II matrices arise naturally in connection with some interesting combinatorial and geometric structures.
\end{abstract}

\section{Introduction}

If $M$ and $N$ are matrices of the same order, their \textsl{Schur product}\/ is the matrix 
$M\circ N$, defined by the condition
$$
(M\circ N)_{i,j}=M_{i,j}N_{i,j}.
$$
The Schur product is commutative and associative, with an identity element $J$, the all-ones matrix.  If $M\circ N=J$ we say that $N$ is the \textsl{Schur inverse} of $M$, and denote it $M\snv$.

A \textsl{type-II matrix}\/ is a Schur invertible $n\times n$ matrix $W$
over $\cx$ such that
$$ 
WW\snt =nI.
$$
This condition implies that $W\inv$ exists and
$$
W\snt =nW\inv.
$$

In \cite{MR89m:57005} Jones showed that certain special type-II matrices could be used to construct
so-called spin models, which could in turn be used to construct interesting invariants
of knots and links (including the Jones polynomial).  The main goal of this paper is to show
that type-II matrices are much more common than might be expected: in particular they
arise in connection with a range of combinatorial and geometric 
structures: symmetric designs, sets of equiangular lines and strongly regular graphs.

\section{The Basics}

We offer some examples of type-II matrices.  First
$$
\pmat{1&1\\ 1&-1}
$$
is a symmetric type-II matrix.  If $\om$ is a primitive cube root of unity then 
$$
\pmat{1&1&\om\\ \om&1&1\\ 1&\om&1}
$$
is also type-II.  For any non-zero complex number $t$, the matrix
$$
W=\pmat{
	1&1&1&1\\
	1&1&-1&-1\\
	1&-1&t&-t\\
	1&-1&-t&t}
$$
is type-II.  Next we have the \textsl{Potts models}: if $W$ is $n\times n$ and
$$
W=(t-1)I+J,
$$
then
\begin{align*}
WW\snt &=((t-1)I+J)((t\inv-1)I+J)\\
	&=((2-t-t\inv)I+(n-2+t+t\inv)J,
\end{align*}
whence it follows that $W$ is type-II whenever $2-t-t\inv=n$, i.e.,
whenever $t$ is a root of the quadratic
$$
t^2+(n-2)t+1.
$$

As the first example suggests, any Hadamard matrix is a type-II matrix, and it is not 
unreasonable to view type-II matrices as a generalisation of Hadamard matrices.

The Kronecker product of two type-II matrices is a type-II matrix; this provides another easy 
way to increase the supply of examples.  Recall that a \textsl{monomial matrix} is the product 
of a permutation matrix and a diagonal matrix.  It is straighforward to verify that if $W$ is 
type-II and $M$ and $N$ are invertible monomial matrices, then $MWN$ is type-II.  We 
say $W'$ is  \textsl{equivalent} to $W$ if $W'=MWN$, where $M$ and $N$ are invertible monomial matrices.

The transpose $W^T$ is also type-II, as is $W\snv$, but these matrices
may not be equivalent to $W$.  It would be a useful exercise to prove that any $2\times2$ type-II matrix is equivalent to the first example above, any $3\times3$ type-II matrix is equivalent to the second, and any $4\times4$ type-II matrix is equivalent to a matrix in the third family.

\medbreak
Let $W$ be a Schur-invertible matrix, with rows and columns indexed by the set $\Om$, where $|\Om|=n$.
Let the vectors
\[
e_a,\ a\in\Om
\]
denote the standard basis for $\cx^\Om$.   We define a set of $n^2$ vectors in $\cx^n$ as follows.
\[
\yab:=We_a\circ W\snv e_b.
\]
We can view $\yab$ as the Schur ratio of the $a$- and $b$-columns of $W$.
The \textsl{Nomura algebra} $\nomw$ of $W$ consists of the set of $n\times n$ complex matrices $M$ such that each of the $n^2$ vectors $Y_{a,b}$ is an eigenvector for $M$.  The Nomura algebra is non-empty, because it always contains $I$.  

\begin{lemma}
\label{lem:JinNom}
Let $W$ be a Schur invertible and invertible matrix.  Then $W$ is a type-II matrix if and only if 
$J\in\nomw$.
\end{lemma}

\proof
Let $D_a$ be the $n\times n$ diagonal matrix such that
\[
(D_a)_{i,i} :=W_{i,a}.
\]
Since $W$ is invertible, its columns
\[
We_b,\quad b\in\Om
\]
are linearly independent.  Since $D_a$ is invertible and
\[
\yab =D_a\inv We_b,
\]
we see that the vectors
\[
\yab,\quad b\in\Om
\]
are linearly independent and consequently they form a basis for $\cx^n$.  

Now $\yy aa= {\bf 1}$, so $J\in \nomw$ if and only if 
\[
J\yab = n\delta_{ab}\yy ab,
\]
equivalently
\[
\sum_{r} \frac{W_{r,a}}{W_{r,b}} = (W\snt W)_{b,a} = n\delta_{b,a}
\]
for all $a, b \in \Om$. 
\qed

It follows that if $W$ is a type-II matrix of order $n\times n$, then $\nomw$ contains $I$ and $J$
and $\dim\nomw\ge2$ when $n\ge2$.  We say that $\nomw$ is \textsl{trivial} if $\dim\nomw=2$.
All the work in this paper is motivated by the desire to find type-II matrices with non-trivial Nomura algebras.  One reason is that if $W$ is a type-II matrix and $W\in\nomw$, then
we may use $W$ to construct a link invariant.  The Potts model, which we mentioned above,
has this property and the corresponding link invariants are evaluations of the Jones polynomial.
For more on this connection, see \cite{MR99f:05125} and \cite{MR98a:57010}.

A type-II matrix $W$ such that $W\in\nomw$ is known as a \textsl{spin model}.  The Potts model
aside, very few interesting spin models are known.  If $W$ is a spin model other than
the Potts model, then $\nomw$ contains $I$, $J$ and $W$, and therefore $\dim\nomw\ge3$.  Spin models have proved very difficult to find.  Hence we are lead to search for type-II matrices whose Nomura algebras are non-trivial.
For reasons that are not at all clear, even these seem to be scarce.

The previous discussion glosses over one point.  If $W_1$ and $W_2$ are type-II matrices,
then the Nomura algebra of $W_1\otimes W_2$ is the tensor product of the Nomura algebras
of $W_1$ and $W_2$.  Since
\[
\dim(\cN_{W_1\otimes W_2}) =\dim(\cN_{W_1})\dim(\cN_{W_2}),
\]
the Nomura algebra of $W_1\otimes W_2$ is always non-trivial.  However the corresponding
link invariants are of no interest, since they are built in an obvious way from the invariants
belonging to the factors.  Therefore our search is actually for type-II matrices which have non-trivial Nomura algebra and which are not equivalent to Kronecker products of type-II matrices.

\section{Nomura Algebras}

We have introduced type-II matrices and their Nomura algebras.  Now we describe
the connection between type-II matrices and combinatorics; the connection is mediated by
association schemes.

Let $W$ be a type-II matrix or order $n\times n$.
We saw in the previous section that
\[
\yab,\quad b\in\Om
\]
form a basis for $\cx^n$.  
If $M\in\nomw$,
then the matrix representing $M$ relative to this basis is diagonal, from which we conclude
that if $M,N\in\nomw$ then $MN=NM$.  In other words, the Nomura algebra of a type-II
matrix is commutative.
We will also see that it is closed under the Schur product.

Let $W$ be a type-II matrix, with rows and columns indexed by the set $\Om$, where $|\Om|=n$.  
If $M\in\nomw$, there is an $n\times n$ matrix $\thmw$ such that
\[
M\yab =(\thmw)_{a,b}\yab.
\]
We call $\thmw$ the \textsl{matrix of eigenvalues} of $M$.  (When no confusion will result, we 
write $\thm$ rather than $\thmw$.)  Note that
\[
\Th(MN) =\thm\circ\Th(N).
\]
Also $\Th$ is an injective linear map from $\nomw$ into the space of $n\times n$ complex matrices.

We define a second family of $n^2$ of vectors in $\cx^n$ as follows.
$$
\ypab :=W^Te_a\circ W^{(-)T}e_b
$$
Thus $\ypab$ is the Schur ratio of two columns of $W^T$, and so the set of matrices with the vectors $\ypab$ as eigenvectors is $\nom{W^T}$.  The following critical result is due to Nomura \cite{MR98g:05158}; it shows that the image of $\nomw$ under $\Th$ is contained in
$\nom{W^T}$.

\begin{theorem}
\label{thm:nomura}
If $M\in\nomw$ then
$$
\thm Y'_{s,r} = n M_{r,s} Y'_{s,r}.
$$
\end{theorem}

\proof
Suppose
$$
F_i := \frac1n Y_{u,i}Y_{i,u}^T.
$$
We verify easily that
$$
F_i F_j =\de_{i,j} F_i,
$$
which shows that the $F_i$'s form an orthogonal set of $n$ idempotents.  We note that $\rk(F_i)=1$ and $\tr(F_i)=1$.  As the $F_i$'s commute it follows that $\sum_i F_i$ is an idempotent matrix with trace equal to $n$; hence
$$
\sum_i F_i =I.
$$
We have
$$
MF_i ={1\over n} MY_{u,i}Y_{i,u}^T =(\thm)_{u,i} F_i.
$$
Summing this over $i$ in $\Om$, recalling that $\sum_i F_i=I$, we get
\begin{equation}
\label{mfi}
M = \sum_i (\thm)_{u,i} F_i.
\end{equation}

Now
$$
(F_i)_{r,s}={1\over n}{W_{r,u}\over W_{r,i}}{W_{s,i}\over W_{s,u}}
	={1\over n}{W_{r,u}\over W_{s,u}}{W_{s,i}\over W_{r,i}}
$$
and therefore, by (\ref{mfi}),
$$
M_{r,s}={1\over n}{W_{r,u}\over W_{s,u}}
	\sum_i (\thm)_{u,i}{W_{s,i}\over W_{r,i}}.
$$
Hence
$$
nM_{r,s}(Y'_{s,r})_u = (\thm Y'_{s,r})_u,
$$
which implies the theorem.\qed
It is an easy consequence that the $\nomw$ is
closed under the Schur product.

We describe a simple way to test if two eigenvectors $\yy ab$'s 
belong to the same eigenspace of $\nomw$.
\begin{lemma}
\label{lem:overlap}
If $Y^T_{a,u}Y_{b,c}\neq 0$ then $(\thm)_{u,a}=(\thm)_{b,c}$. 
\end{lemma}

\proof
It follows from $\sum_i F_i = I$ that
$$
Y_{b,c}= \frac{1}{n}\sum_i (Y^T_{i,u}Y_{b,c})Y_{u,i}.
$$
So
$$
(\thm)_{b,c}Y_{b,c} =MY_{b,c}
	={1\over n}\sum_i (Y^T_{i,u}Y_{b,c}) (\thm)_{u,i} Y_{u,i}.
$$
Multiply both sides of this by $Y^T_{a,u}$ to get
\begin{align*}
(\thm)_{b,c}Y^T_{a,u}Y_{b,c}
	&={1\over n}(Y^T_{a,u}Y_{b,c}) (\thm)_{u,a}Y^T_{a,u}Y_{u,a}\\
	&=Y^T_{a,u}Y_{b,c} (\thm)_{u,a}.
\end{align*}
If $Y^T_{a,u}Y_{b,c}\ne0$, this implies that
$(\thm)_{u,a}=(\thm)_{b,c}$.\qed

\section{Association Schemes}
We recall some definitions.  An \textsl{association scheme} with $d$ classes is a collection 
$\cA$ of $01$-matrices  $\seq A0d$ of order $n\times n$ such that:
\begin{enumerate}[(a)]
\item
$A_0=I$.
\item
$\sum_i A_i =J$.
\item
$A_i^T\in\cA$ for $i=0,\ldots,d$.
\item
The product $A_iA_j$ lies in the span of $\cA$, for all $i$ and $j$.
\item
$A_iA_j=A_jA_i$.
\end{enumerate}
The matrices $A_i$ are the adjacency matrices of directed graphs whose arc 
sets partition the arcs of the complete directed graph on $n$ vertices.  
It follows from the axioms that $A_iJ=JA_i$, whence each directed graph 
is regular.  The span of $\cA$ is called the \textsl{Bose-Mesner algebra} 
of the association scheme.  Since the $A_i$ are $01$-matrices and sum to $J$, 
they form a basis for $\cA$; since the set consisting of $\cA$ and 
the zero matrix is closed under Schur product, 
it follows that the span of $\cA$ is closed under the Schur product.  
The axioms also insure that the Bose-Mesner algebra is closed under 
transpose and under matrix multiplication.
On the other hand, a vector space of matrices is the Bose-Mesner algebra of an
association scheme if it contains $I$ and $J$ and is closed under transpose, 
matrix and Schur product, and it is commutative 
with respect to matrix multiplication.  See \cite{MR1002568} for details.

The simplest example of an association scheme arises if we take $d=1$ and $A_1=J-I$.  This is the
association scheme of the complete graph (and the Nomura algebra of a Potts model).

\begin{corollary}
If $W$ is a type-II matrix, then $\nomw$ is the Bose-Mesner algebra of an association scheme.
\end{corollary}

\proof
It is immediate from its definition that $\nomw$ is closed under 
and is commutative with respect to matrix multiplication.
By \lref{JinNom}, $\nomw$ contains $I$ and $J$.
We show that it is also closed under transpose and Schur multiplication. 

\tref{nomura} yields that if $M\in\nomw$ then
\[
\Th_{W^T}(\Th_W(M)) =nM^T.
\]
Since $\Th_W(M)\in\nom{W^T}$ we see that $\Th_{W^T}(\Th_W(M))\in\nomw$.  Therefore $\nomw$ is closed under transpose.  We also see that $\Th_{W^T}$ is surjective.

We saw that if $M,N\in\nom{W^T}$ then
\[
\Th_{W^T}(MN)=\Th_{W^T}(M)\circ\Th_{W^T}(N).
\]
This shows that $\Th_{W^T}(\nom{W^T})$ is closed under Schur multiplication.  Since 
$\Th_{W^T}$ is surjective, it follows that $\nomw$ is closed under Schur multiplication.\qed

If $W$ is a type-II matrix with algebra $\nomw$ then, as noted before, 
$W$ determines a spin model if and only if $W$ lies in $\nomw$. 
As any type-II matrix equivalent to $W$ has the isomorphic Nomura algebra, \cite{MR99f:05125}, 
we may concentrate on the matrices $W$ that lie in their Nomura algebra.  If $W\in\nomw$ then
\begin{enumerate}[(a)]
\item 
$W$ is normal.
\item
The diagonal of $W$ is constant, that is, $W\circ I=cI$ for some $c$.
\item
The row and column sums of $W$ are all equal.
\end{enumerate}
These conditions hold because they are satisfied by any matrix in a Bose-Mesner algebra.

One consequence of Nomura's theorem is that, when searching for spin models, we can restrict 
ourselves to type-II matrices that lie in the Bose-Mesner algebra of an association scheme.  This is 
important because there may be uncountably many type-II matrices of a given order $n$, but there 
are only finitely many association schemes of order $n$.  Hence our search space is considerably
restricted.

\section{Hadamard Matrices}

A Hadamard matrix is a $\pm1$-matrix of order $n\times n$ such that 
\[
H^TH=nI.
\]
Since $H\circ H=J$ it follows that $H$ is a type-II matrix.  Hadamard matrices have long been
of interest to combinatorialists.  Since they are the simplest examples of type-II matrices,
we summarise what is known about their Nomura algebras here.

\begin{lemma}
If $W$ is real then all matrices in $\nomw$ are symmetric.
\end{lemma}

\proof
If $W$ is real then the eigenvectors $Y_{a,b}$ are real.  Hence the Schur idempotents of the 
scheme have only real eigenvalues.  Since $\nomw$ is closed under transposes and is a commutative  algebra, the Schur idempotents are real normal matrices.  A real normal matrix is symmetric if and only if its eigenvalues are real.\qed

The following is a new proof of a result due to Jaeger et al 
\cite{MR99f:05125}.

\begin{lemma}
\label{had8}
Let $W$ be a Hadamard matrix of order $n$.  If $\nomw$ is non-trivial, then $n$ is divisible by eight.
\end{lemma}

\proof
Let $w_i$ denote $We_i$.  Normalise $W$ so that $w_1=\one$ and assume
$1$, $i$, $j$ and $k$ are distinct.  Then
$$
(w_1+w_i)\circ(w_1+w_j)\circ(w_1+w_k)
$$
is the Schur product of three vectors with entries $0,\pm2$.  The sum
of the entries of this vector is 
\begin{multline}
\ip\one{w_1^{\circ 3}} +\ip\one{w_1^{\circ2}\circ(w_i+w_j+w_k)}\hfill\\
		\hfill+\ip\one{w_1\circ(w_i\circ w_j+w_i\circ w_k+w_j\circ w_k)}
		+\ip\one{w_i\circ w_j\circ w_k}
\end{multline}
Since $W$ is a Hadamard matrix, the second and third terms here are
zero, whence we deduce that, modulo 8,
$$
n+\ip\one{w_i\circ w_j\circ w_k}=0
$$
and therefore, if $n$ is not divisible by 8, then $\yy i1=w_i$ cannot be
orthogonal to $\yy jk =w_j\circ w_k$.\qed

If $H$ is a Hadamard matrix of order less than 32, its Nomura algebra is a product of Potts models. (Unpublished computations by Allan Roberts and the second author.)

Hadamard matrices form a special class of a more general class of type-II matrices.
A complex matrix is \textsl{flat} if all its entries have the same absolute value.
The following result is easy to prove.

\begin{lemma}
For an $n\times n$ matrix, any two of the following statements imply
the third:
\begin{enumerate}[(a)]
\item $W$ is a type-II matrix.
\item $n^{-1/2}W$ is unitary.
\item $|W_{i,j}|=1$ for all $i$ and $j$.\qed
\end{enumerate}
\end{lemma}

In other words, a unitary matrix is type-II if and only if it is flat.  The character table of
an abelian group is flat, type-II and unitary.  Flat unitary matrices appear in quantum physics
in connection to mutually unbiased sets of orthogonal bases.

\section{Symmetric Designs}

%In the previous section we considered type-II matrices with exactly two distinct eigenvalues.  
%Now we impose a constraint on the entries, rather than the eigenvalues.
We consider type-II matrices with exactly two distinct entries,
that are not Hadamard matrices.

\begin{theorem}
\label{thm:symd}
Suppose $W=aJ+(b-a)N$, where $N$ is a $01$-matrix 
and $a\neq \pm b$.  
Then $W$ is type II if and only if $N$ is the incidence matrix of a symmetric design.
\end{theorem}

\proof
Let $N$ be the incidence matrix of a symmetric $(v,k,\la)$-design, and
let $W$ be given by
$$
W =J+(t-1)N,
$$
where
$$
t ={1\over2(k-\la)}\left(2(k-\la)-v\pm\sqrt{v(v-4(k-\la))}\right).
$$
We show that $W$ is a type-II matrix.

We have
$$
W\snv =(t\inv-1)N+J
$$
and, as $NJ=N^TJ=kJ$ and $J^2=vJ$,
\begin{align*}
WW\snt &=(t-1)(t\inv-1)NN^T +(k(t+t\inv-2)+v)J\\
	&=(t-1)(t\inv-1)(k-\la)I+((k-\la)(t+t\inv-2)+v)J.
\end{align*}
The coefficient of $J$ is zero if
$$
(k-\la)(t-1)^2+v(t-1)+v=0,
$$
which yields sufficiency.

We now prove the converse.  If $W$ has exactly two distinct entries,
there is no harm in assuming that we have
$$
W=J+(t-1)N
$$
for some $01$-matrix $N$ and some complex number $t$ such that
$t\ne \pm 1$.  Then $W\snt=J+(t\inv-1)N^T$ and so, if $W$ is $v\times v$, we
have
$$
WW\snt =vJ+(t-1)NJ+(t\inv-1)JN^T+(t-1)(t\inv-1)NN^T.
$$
Since $WW\snt=vI$ and $NN^T$ is symmetric, this implies that
$$
M := (t-1)NJ+(t\inv-1)JN^T
$$
is symmetric.  We work with this.  Note that this equation yields
$$
M-M^T =(t-t\inv)NJ+(t\inv-t)JN^T =(t-t\inv)(NJ-(NJ)^T).
$$
Since $M=M^T$ and $t\ne \pm 1$, this forces us to conclude that $NJ$ is
symmetric.  Hence there is a positive integer $k$ such that
$$
NJ=JN^T=kJ.
$$

Returning to our expression for $WW\snt$, we now have
\begin{equation}
WW\snt =(v+k(t+t\inv-2))J+(2-t-t\inv)NN^T.
\label{eqnSymm}
\end{equation}
Since $(2-t-t\inv)=-(t-1)^2/t$ and $t\ne1$, it follows that $NN^T$ is
a linear combination of $I$ and $J$, and consequently $N$ is the
incidence matrix of a symmetric design.
\qed

Note that if $v+k(t+t\inv-2)=0$ in \ref{eqnSymm})
then we get $NN^T=kI$.  
Since $N$ is a square $01$-matrix, $NN^T=kI$ only when $k=1$.
In this case, $N$ 
is the incidence matrix of the complement of the complete design,
and $W=J+(t-1)N$ is equivalent to the Potts model.

If $H$ is a Hadamard matrix, we may multiply it fore and aft by
diagonal matrices, thus setting all entries in the first row and
column to 1.  If $H_1$ is the matrix we get from this by deleting the
first row and column, then
$$
{1\over2}(H_1+J)
$$
is the incidence matrix of a symmetric design.  This gives a large
class of examples of symmetric designs.

\begin{lemma}
\label{lem:symm}
Suppose $W$ is a type-II matrix of the form $(t-1)N+J$, where $N$ is
the incidence matrix of a symmetric $(v,k,\la)$-design.  If $v>3$,
then all matrices in $\nom{W^T}$ are symmetric.
\end{lemma}

\proof
We show that $\ip{Y_{i,j}}{Y_{i,j}}\ne0$ when $v>3$.  By
\lref{overlap}, it follows that $\Th(M)_{i,j}=\Th(M)_{j,i}$ for all $M$
in $\nomw$ and for all $i$ and $j$.

We have
\begin{align*}
\ip{Y_{i,j}}{Y_{i,j}} &=(k-\la)(t^2+t^{-2})+v-2k+2\la\\
	&=(k-\la)(t^2-2+t^{-2})+v,
\end{align*}
and so, if $\ip{Y_{i,j}}{Y_{i,j}}=0$ then
$$
t^2-2+t^{-2} ={-v\over(k-\la)}.
$$
From our computations in the proof of the previous theorem,
\begin{equation}
\label{klatv}
(k-\la)(t-2+t\inv)+v=0,
\end{equation}
and so
$$
t-2+t\inv ={-v\over(k-\la)}.
$$
As 
$$
t^2-2+t^{-2}  =(t-2+t\inv)(t+2+t\inv),
$$
these equations imply that, if $\ip{Y_{i,j}}{Y_{i,j}}=0$, then
$$
t+1+t\inv =0,
$$
whence \eqref{klatv} implies that $v=3(k-\la)$.  

Since $v(v-1)\la=vk(k-1)$, if $v=3(k-\la)$, then
$$
k^2=k+(v-1)\la=(3\la+1)(k-\la)
$$
and therefore
$$
k^2-(3\la+1)k+3\la^2+\la=0
$$
This discriminant of this quadratic is
$$
1+2\la-3\la^2 =(1-\la)(1+3\la),
$$
which is negative if $\la>1$.  The lemma follows.\qed

\begin{lemma}
Let $N$ be the incidence matrix of a symmetric design, and let $W$ be
a type-II matrix of the form $(t-1)N+J$.  If $t\ne-1$, then the
difference of two distinct columns of $N$ is an eigenvector for the
Nomura algebra of $W$.
\end{lemma}

\proof
If $u$ is a point in the design and $\al$ and $\be $ are the $i$-th
and $j$-th blocks in the design, then
$$
(\yij)_u=\begin{cases}
	t,& \mathrm{if}\ u\in\al\diff\be;\\
	t\inv,&\mathrm{if}\ u\in\be\diff\al;\\
	1, &\mathrm{otherwise}.
	\end{cases}
$$
By the previous lemma, $\yij$ and $Y_{j,i}$ have the same eigenvalues
for any matrix in $\nomw$.  Therefore the vector
$$
(t-t\inv)\inv(\yij-Y_{j,i})
$$
is an eigenvector for each matrix in $\nomw$, but this vector is just
the difference of the $i$-th and $j$-th columns of $N$.\qed

We note that if $t=-1$ then $(t-1)N+J$ is type II if and only if it is a
Hadamard matrix.
The previous lemmas
lead to the following disappointing consequence.

\begin{theorem}
Suppose $W$ is a type-II matrix of the form $(t-1)N+J$, where $N$ is
the incidence matrix of a symmetric $(v,k,\la)$-design.  
If $v>3$ and $t\ne-1$, then the Nomura algebra of $W$ is trivial.
\end{theorem}

Let $Z_{i,j}:=Ne_i-Ne_j$ for some $i \neq j$.  
If $k$ is distinct from $i$ and $j$ then
$$
\ip{Z_{i,j}}{Ne_k} =\ip{Ne_i}{Ne_k}-\ip{Ne_j}{Ne_k}=\la-\la=0
$$
while
$$
\ip{Z_{i,j}}{Ne_i} =k-\la.
$$
We conclude that $\ip{Z_{i,j}}{Z_{i,k}}=k-\la$ and therefore at least one of
$$
Y^T_{i,k}\yij, \quad Y^T_{k,i}\yij,\quad
Y^T_{i,k}Y_{j,i} \quad \text{and}\quad Y^T_{k,i}Y_{j,i}
$$
is non-zero.  It follows from \lref{overlap} and \lref{symm} that
$$
\Th(M)_{i,k}=\Th(M)_{i,j}
$$ 
for any matrix $M$ from $\nomw$.  It follows that $\nomw$ must be trivial.\qed

\section{Equiangular Lines}

We consider sets of lines in $\cx^d$.  A set of lines in $\cx^d$ spanned by the unit vectors 
$\seq x1n$ is \textsl{equiangular} if there is a real number $\al$ such that
\[
|\ip{x_i}{x_j}|=\al
\]
whenever $i\ne j$.  Note that it is reasonable to take the absolute value here, because if $\la\in\cx$ and $|\la|=1$ then $\la x_i$ and $x_i$ are unit vectors spanning the same line.
We will refer to $\al$ as the \textsl{angle} between the lines.  We are also interested in 
equiangular sets of lines in $\re^d$; the above definition still works in this case.  We have the following result, due to \cite{MR92m:01098}.

\begin{theorem}
If there is a set of $n$ equiangular lines in $\cx^d$ or $\re^d$ with angle $\al$ and 
$d\al^2<1$, then
\[
n\le\frac{d(1-\al^2)}{1-d\al^2}.
\]
\end{theorem}

\proof
Suppose $\seq x1n$ are unit vectors spanning a set of equiangular lines in $\cx^d$
and suppose $X_i :=x_ix_i^*$.  Then $X_i$ is a Hermitian matrix that represents
orthogonal projection onto the line spanned by $x_i$.  Assume that $|\ip{x_i}{x_j}|=\al$
when $i\ne j$.  The space of Hermitian matrices is a real inner product space with inner 
product $\ip{X}{Y}$ given by
\[
\ip{X}{Y} =\tr(XY).
\]
Then $\ip{X_i}{X_i}=1$ and if $i\ne j$ then
\begin{align*}
\ip{X_i}{X_j} =\tr(X_iX_j) &=\tr(x_ix_i^*x_jx_j^*)\\
	&=\tr(x_j^*x_ix_i^*x_j)\\
	&=|x_i^*x_j|^2\\
	&=\al^2.
\end{align*}

If 
\[
Z :=\sum_i X_i
\]
then
\[
\ip{Z}{Z} =n+(n^2-n)\al^2
\]
and if $\ga\in\re$, then
\[
\ip{Z-\ga I}{Z-\ga I}  =n+(n^2-n)\al^2 -2\ga n +\ga^2d.
\]
Here the right side is a quadratic in $\ga$, and is non-negative for all real $\ga$.
Its minimum value occurs when $\ga =n/d$, which implies that
\[
-\frac{n^2}{d} +n(1+\al^2(n-1))\ge 0.
\]
The theorem follows from this.\qed

Note that the above proof still works if we replace $\cx$ by $\re$ and `Hermitian' by
`symmetric'.

We say a set of lines is \textsl{tight} if equality holds in the bound of the previous theorem.
We say that an $n\times n$ matrix $C$ is a \textsl{generalized conference matrix} if:
\begin{enumerate}[(a)]
\item
$C$ is Hermitian
\item
$C_{i,i}=0$ for all $i$.
\item
$|C_{i,j}|=1$ if $i\ne j$.
\item
The minimal polynomial of $C$ is quadratic.
\end{enumerate}
Note that a conference matrix is an $n\times n$ matrix with diagonal entries zero and
off-diagonal entries $\pm1$, such that $CC^T=(n-1)I$.  It is known that a conference matrix
is equivalent to a symmetric or skew symmetric conference matrix.  If $C$ is symmetric
then it is Hermitian and $C^2-(n-1)I=0$.  If $C$ is skew symmetric, then $iC$ is
Hermitian and $(iC)^2-(n-1)I=0$.

\begin{corollary}
Suppose $\seq x1n$ are unit vectors that span a set of equiangular lines in $\cx^d$ with 
angle $\al$ and Gram matrix $G$, and suppose $G=I+\al C$.
Then the set of lines is tight if and only if $C$ is a generalized conference matrix.
\end{corollary}

\proof
Suppose $\seq x1n$ span a set of equiangular lines in $\cx^d$, let $X_i$ be the orthogonal projection onto the line spanned by $x_i$ and set $Z=\sum_i X_i$.  
If this set of lines is tight, then
\[
\ip{Z-\ga I}{Z-\ga I}  =0
\]
and consequently
\[
\sum_i X_i =\frac{n}{d}I.
\]
Let $U$ be the $n\times d$ matrix with $i$-th row equal to $x_i^*$.  Then
\[
U^*U =\sum_i X_i =\frac{n}{d}I.
\]
Now $G :=UU^*$ is the Gram matrix of the unit vectors $\seq x1n$; since
$UU^*$ and $U^*U$ have the same non-zero eigenvalues with the same multiplicities
it follows that the eigenvalues of $G$ are 0 and $n/d$.  Since our set of lines is equiangular,
we may write
\[
G =I+\al C.
\]
Here $C$ is Hermitian, its diagonal entries are zero, its off-diagonal entries all have 
absolute value 1, and its minimal polynomial is quadratic.  Thus it is a generalized
conference matrix.

For the converse, suppose that $C$ is a non-zero Hermitian matrix with zero diagonal and
\[
C^2-\be C-\ga I =0.
\]
Then the diagonal entries of $C^2$ are positive, whence $\ga\ne0$ and $C$ is invertible.
If $\tau$ is the least eigenvalue of $C$, then
\[
G := I-\frac1\tau C
\]
is Hermitian and all its eigenvalues non-negative.  Assume $\rk(G)=d$.  Since $\tr(G)=n$ it follows that the eigenvalues of $G$ are 0 and $n/d$.  Hence there is an $n\times d$ matrix $U$ such that
\[
U^*U =\frac{n}{d}I,\quad UU^* =G.
\]
Thus $G$ is Gram matrix of the columns of $U^*$, and so these columns span a set of
equiangular lines in $\cx^d$.  Since $U^*U=(n/d)I$, the set of lines is tight.\qed

Conditions (a) and (c) in the definition of generalized conference matrix imply that $(C^2)_{i,i}=(n-1)I$, whence the minimal polynomial of $C$ has the form 
$z^2-\be z-(n-1)$, for some $\be$.

\begin{theorem}
Suppose $C$ is a generalized conference matrix of order $n\times n$ with minimal
polynomial $z^2-\be z -(n-1)$.  If $t+t\inv+\be=0$, then $tI+C$ is type II.
\end{theorem}

\proof
If $C$ is a generalized conference matrix, then 
\[
(tI+C)\snt =t\inv I+C
\]
and therefore
\begin{align*}
(tI+C)(tI+C)^{(-)T} &=I+t\inv C +tC\snt+C C\snt\\
	&=I+(t+t\inv)C+C^2\\
	&=I+(t+t\inv)C+\be I +(n-1)I\\
	&=nI+(t+t\inv+\be)C.
\end{align*}
Hence $tI+C$ is type-II if
\[
t+t\inv+\be =0.\qed
\]

We derive a converse to this result, under weaker conditions.

\begin{theorem}
Let $W$ be a type-II matrix with all diagonal entries equal to $c$ and with quadratic minimal polynomial.  If $W-cI$ is Hermitian, it is a scalar multiple of a generalized conference matrix.
\end{theorem}

\proof
Suppose that $W$ is $n\times n$ and
$$
W^2-\be W -\ga I=0.
$$
Since $W$ is invertible, $\ga\ne0$ and 
$$
W\inv  =-\frac{1}{\ga}(\be I-W).
$$
Hence
$$
J =nW\circ W^{-T} = -\frac{n}{\ga}(\be W\circ I-W\circ W^T),
$$
from which we find that
\begin{equation}
\label{aaaij}
W\circ W^T =\be W\circ I +\frac{\ga}{ n}J.
\end{equation}
It follows that all off-diagonal entries of $W$ have the same absolute value (namely 
$\sqrt{\ga/n}$).\qed

\section{Strongly Regular Graphs}

A graph $X$ is strongly regular if it is not complete and there are integers $k$, $a$ and $c$
such that the number of common neighbours of an ordered pair of vertices $(u,v)$ is
$k$, $a$ or $c$ according as $u$ and $v$ are equal, adjacent or distinct and not
adjacent.  Trivial examples are provided by the graphs $mK_n$ and their complements.
The Petersen graph provides a less trivial example.  A strongly regular graph $X$ is
\textsl{primitive} if both $X$ and its complement are connected; an imprimitive
strongly regular graph is isomorphic to $mK_n$ or its complement.
A strongly regular graph $X$ gives rise to an
association scheme with two classes, corresponding to $X$ and its complement.  Conversely
each association scheme with two classes determines a complementary pair of strongly regular
graphs.  

\begin{theorem}
Let $X$ be a primitive strongly regular graph with $v$ vertices, valency $k$,
and eigenvalues $k$, $\th$ and $\tau$,
where $\th>\tau$.  Let $A_1$ be the adjacency matrix of $X$ and $A_2$ the 
adjacency matrix of its complement.  Suppose
\[
W := I+xA_1+yA_2.
\]
Then $W$ is a type-II matrix if and only if one of the following holds
\begin{enumerate}[(a)]
\item
$y=x ={1\over2}(2-v\pm\sqrt{v^2-4v})$.
\item
$x=1$ and 
$y=1+\frac{1}{2(\bar{k}-\la)}(-v\pm\sqrt{v^2-4(\bar{k}-\la)v})$ and $A_2$ is the incidence 
matrix of a symmetric $(v, \bar{k}, \la)$-design where $\bar{k} = v-k-1$.
\item
$x=-1$ and $y =\frac12(\la\pm\sqrt{\la^2-4})$ (where
$\la =(1+\th \tau)\inv (2-2\th \tau-v)$), 
and $A_1$ is the incidence matrix of a symmetric design.
\item
$x+x\inv$ is a zero of the quadratic $z^2-\al z+\be-2$
with 
 \begin{align*}
\al &={1\over \th \tau}[v(\th +\tau+1)+(\th +\tau)^2],\\
\be &={1\over \th \tau}[-v-v(1+\th +\tau)^2+2\th ^2+2\th \tau+2\tau^2]
\end{align*}
and
\[
y=\frac{1}{(x-x\inv)}\left(\frac{\th \tau x-1}{(\th +1)(\tau +1)}(x+x\inv-2+v)-(v-2)x-2\right).
\]
\end{enumerate}
\end{theorem}

\proof
We use $\ell$ to denote valency $v-1-k$ of the complement of $X$.
Then the eigenvalues of $A_2$ are $v-1-k$, $-1-\tau $ and $-1-\th $
and the equation $WW\snt=vI$ is equivalent to 
\begin{align*}
(1+kx+\ell y)(1+kx\inv+\ell y\inv)&=v,\\
(1+\th x+(-\th -1)y)(1+\th x\inv+(-\th -1)y\inv)&=v,\\
(1+\tau x+(-\tau -1)y)(1+\tau x\inv+(-\tau -1)y\inv)&=v. 
\end{align*}
Note that this set of equations is invariant under the substitutions
$$
x\mapsto x\inv,\quad y\mapsto y\inv
$$
and also under the substitutions
$$
x\mapsto y,\quad y\mapsto x,\quad \th \mapsto-\th -1,\quad \tau \mapsto-\tau -1.
$$

The missing details in the following calculations were performed in
Maple. 

If we set
$$
X := x+{1\over x},\quad Y :=y+{1\over y},
	\quad Z:= {x\over y}+{y\over x}
$$
then, from our three equations we get
\begin{align}
k\ell Z +kX +\ell Y&=v-1-k^2-\ell^2,\nonumber\\
-\th (\th +1)Z +\th X -(\th +1)Y&=v-1-\th ^2-(\th +1)^2,\label{thxy}\\
-\tau (\tau +1)Z +\tau X -(\tau +1)Y&=v-1-\tau ^2-(\tau +1)^2.\label{tauxy}
\end{align}
These three equations are linearly dependent: if $\th $ has multiplicity
$m$ and $\tau $ has multiplicity $n$ as an eigenvalue of $A_1$, then the
first equation plus $m$ times the second plus $n$ times the third is
zero.  In fact, our three equations are equivalent to the following
pair.
\begin{align}
Y-2+v&={\th \tau \over(\th +1)(\tau +1)}(X-2+v),\label{yvxv}\\
Z-2&={1\over(\th +1)(\tau +1)}(X-2+v).\label{zxv}
\end{align}
Given the definitions of $Y$ and $Z$, we can view this as a pair of
linear equations in $y$ and $y\inv$, whence we find that
$$
y(x-x\inv)={\th \tau x-1\over(\th +1)(\tau +1)}(x+x\inv-2+v)-(v-2)x-2.
$$
Assume $x^2\ne1$.  
If we define
\[
p(x) := \tau \th x^3+(1-v+2 \th+2 \tau-\th v-\tau v) x^2 -
(2 \th+\tau \th+2 \tau+v) x-1,
\]
then \eqref{yvxv} and \eqref{zxv} hold if and only if
$$
y={p(x)\over(\th +1)(\tau +1)(x^2-1)},\qquad 
	y\inv={-x^2p(x\inv)\over(\th +1)(\tau +1)(x^2-1)}.
$$
 
Then the previous expressions for $y$ and $y\inv$ hold if and only if
\[
-x^2p(x)p(x\inv)=[(\th+1)(\tau+1)(x^2-1)]^2.
\]
We deduce that $x$ must be a root of the polynomial
\begin{equation}
(x^2+(v-2)x+1)(x^4-\al x^3 +\be x^2-\al x+1)\label{hex}
\end{equation}
where
\begin{align*}
\al &={1\over \th \tau }[v(\th +\tau +1)+(\th +\tau )^2],\\
\be &={1\over \th \tau }[-v-v(1+\th +\tau )^2+2\th ^2+2\th \tau +2\tau ^2].
\end{align*}

If $x$ is a root of the quadratic factor in \eqref{hex}, then
$X-2+v=0$ and so Equations \eqref{yvxv} and \eqref{zxv} imply that
$Y=2-v$ and $Z=2$.  Since 
$$
Z-2={(x-y)^2\over xy},
$$
it follows that
$$
y=x ={1\over2}(2-v\pm\sqrt{v^2-4v}).
$$
This is the Potts model solution.

We turn to the quartic factor in \eqref{hex}, which is equal to
\[
x^2\left(\left(x+x\inv\right)^2 -\al\left(x+x\inv\right) +\be-2\right).
\]
From this we see that $X$ must  be a zero of the quadratic
\begin{equation}
z^2-\al z+\be-2
\label{qua}
\end{equation}
and thus (d) holds.

To complete the proof we consider the cases where $x^2=1$.  If $x=1$
then \tref{symd} yields that $A_2$ is the incidence matrix of a symmetric design.
So we assume $x=-1$.

Equations \eqref{yvxv} and \eqref{zxv} imply that
\[
Y-2+v =\th \tau (Z-2).
\]
Since $Z=-Y$ if $x=-1$, we find that
\[
(1+\th \tau )Y =2-2\th \tau -v
\]
whence 
\[
y =\frac12(\la\pm\sqrt{\la^2-4}),
\]
where
\[
\la =\frac{2-2\th \tau -v}{1+\th \tau }.
\]
(The denominator cannot be zero because $\tau \le-2$ and $\th \ge 1$ for any 
primitive strongly regular graph.) 

If $x=-1$ then $Z=-Y$ and $X=-2$; if we add equations \eqref{yvxv} and \eqref{zxv}
we get
\[
v-4 =\frac{(\th \tau +1)(v-4)}{(\th +1)(\tau +1)}.
\]
whence we find that $v-4$ or $\th +\tau =0$.  Since, for any strongly regular graph,
\[
A^2-(\th +\tau )A+\th \tau I =(k+\th \tau )J,
\]
we see that if $\th +\tau =0$, then $A^2=-\th \tau I+(k+\th \tau )J$.  Therefore $A$ is the incidence matrix
of a symmetric design (with zero diagonal and symmetric incidence matrix).\qed

Jaeger \cite{MR94e:57008} showed that if $W$ is a spin model then
$X$ is formally self-dual.
If $X$ is formally self-dual then
$v=(\th-\tau)^2$ and
the quadratic (\ref{qua}) becomes
\[
\left(z - {\tau^2-\th^2+2\tau\over \th}\right)
\left(z - {\th^2-\tau^2+2\th\over \tau}\right).
\]
In addition to the Potts model solutions, Equations 
(\ref{yvxv}) and (\ref{zxv}) give
\begin{eqnarray*}
x&=&\frac{1}{2\tau}\left(\th^2-\tau^2+2\th\pm \sqrt{(\th-\tau)(\th-\tau+2)(\th+\tau)(\th+\tau+2)}\right)
\quad \text{and}\\
y&=&\frac{1}{2(\th+1)}\left(\th^2-\tau^2+2(\th+1)\pm
\sqrt{(\th-\tau)(\th-\tau+2)(\th+\tau)(\th+\tau+2)}\right),
\end{eqnarray*}
or
\begin{eqnarray*}
x&=&\frac{1}{2\th}\left(\tau^2-\th^2+2\tau\pm \sqrt{(\th-\tau)(\th-\tau-2)(\th+\tau)(\th+\tau+2)}\right)
\quad \text{and}\\
y&=&\frac{1}{2(\tau+1)}\left(\tau^2-\th^2+2(\tau+1)\pm
\sqrt{(\th-\tau)(\th-\tau-2)(\th+\tau)(\th+\tau+2)}\right).
\end{eqnarray*}
Hence there are at most six type-II matrices, up to equivalence, 
in the Bose-Mesner algebra of a formally self-dual strongly regular graph.

We now determine what happens to the imprimitive strongly regular graphs,
which will arise in the next section.
\begin{theorem}
Let $A_1$ be the adjacency matrix of $mK_{k+1}$ and $A_2$ the 
adjacency matrix of its complement.  
Suppose
\[
W := I+xA_1+yA_2.
\]
Then $W$ is a type-II matrix if and only if one of the following holds
\begin{enumerate}[(a)]
\item
%$y=x ={1\over2}(2-v\pm\sqrt{v^2-4v})$,
$W$ is equivalent to the Potts model,
\item
\[
x=\frac{(kv-2k-1)y^2-(v-2k-2)y-1}{k(1-y^2)}
\]
and
\[
y+y\inv=
\frac{2(k+1)^2-v(k^2+1)}{(k+1)^2-kv}
\]
where $v=m(k+1)$.
\end{enumerate}
\end{theorem}
\proof
The eigenvalues of $A_1$ are $k$ and $-1$, so $\theta=k$ and $\tau=-1$.
The equation $WW\snt=vI$ are equivalent to 
Equations~(\ref{thxy}) and (\ref{tauxy}):
\begin{eqnarray*}
-k(k+1)Z+kX-(k+1)Y &=& v-1-k^2-(k+1)^2\\
X &=& -v+2.
\end{eqnarray*}
Solving this as a pair of linear equations in $x$ and $x\inv$ gives
\begin{equation*}
k(1-y^2)x = (kv-2k-1)y^2-(v-2k-2)y-1. \label{xy}
\end{equation*}
Assume $y^2\neq 1$.
Then Equations~(\ref{thxy}) and (\ref{tauxy}) are equivalent to
\begin{equation*}
x = \frac{p(y)}{k(1-y^2)}
%x = \frac{(kv-2k-1)y^2-(v-2k-2)y-1}{k(1-y^2)}
\end{equation*}
and
\begin{equation*}
x\inv = \frac{-y^2p(y\inv)}{k(1-y^2)}
%\inv = \frac{-y^2 -(v-2k-2)y+ (kv-2k-1)}{k(y^2-1)}.
\end{equation*}
where
\begin{equation*}
p(y)=
(kv-2k-1)y^2-(v-2k-2)y-1.
\end{equation*}
Now these expressions for $x$ and $x\inv$ hold if and only if
\begin{equation*}
-y^2p(y\inv)p(y)=k^2(1-y^2)^2.
\end{equation*}
We deduce that $y$ must be a root of the quartic
\begin{equation*}
\left(y^2+(v-2)y+1\right)
\left(y^2-\be y +1\right)
\end{equation*}
where
\begin{equation*}
\be = \frac{2(k+1)^2-v(k^2+1)}{(k+1)^2-kv}.
\end{equation*}
If $y$ is a root of $y^2+(v-2)y+1$ then we deduce from
Equation~(\ref{xy})  that $x=y$ and $W$ is the Potts model.

If $y=1$ then $Y=2$, $Z=X$ and Equation~(\ref{thxy}) becomes
$X = {-v \over k^2} +2$.  Equations~(\ref{thxy}) and (\ref{tauxy}) imply
$k=1$.  In this case, $A_1$ is a permutation matrix and $W=J + (x-1)A_1$ is 
equivalent to the Potts model.

If $y=-1$ then $Y=-2$, $Z=-X$ and Equation~(\ref{thxy}) becomes
\[
X = {v-2k^2-4k-4 \over k^2+2k}.
\]
Equations~(\ref{thxy}) and (\ref{tauxy}) imply
\[
2-v = {v-2k^2-4k-4 \over k^2+2k},
\]
which leads to $v=4$ and $x=-1$.  In this case, 
$-W=J-2I$ is the Potts model.\qed

\section{Covers of Complete Graphs}

Now we know that the Bose-Mesner of algebra of an association scheme with two
classes contains type-II matrices different from the Potts models.  Given this,
it is natural to ask what happens in schemes with more than two classes;
in this section we consider the next simplest case.  We will see that non-trivial
type-II matrices do arise, and that the amount of effort required to establish this
increases considerably.

We say a graph of diameter $d$ is \textsl{antipodal} if whenever $u$, $v$ and $w$ 
are vertices and
\[
\dist(u,v)=\dist(v,w)=d,
\]
then $u=w$ or $\dist(u,w)=d$.  If $X$ is antipodal, then the relation ``at distance $0$ 
or $d$" is an equivalence relation.  The cube and the line graph of the Petersen graph 
provide two examples with $d=3$.  If $X$ is antipodal with $d=2$, then it is the 
complement of a collection of complete graphs.  If $X$ is an antipodal graph with diameter $d$, 
then its `antipodal classes' form the vertices of a distance-regular graph with the same 
valency and diameter $\lfloor\frac{d}2\rfloor$.

Here we are interested in distance-regular antipodal graphs with diameter three.
To each such graph there is a set of four parameters $(n,r,a_1,c_2)$.  The integer $n$
is the number of antipodal classes, and $r$ is the number of vertices in each class.
If $(u,v)$ is a pair of vertices from $X$ and $\dist(u,v)=1$ then $u$ and $v$ have exactly
$a_1$ common neighbours; if $\dist(u,v)=2$ they have exactly $c_2$
common neighbours.
The value of $a_1$ is determined by $n$, $r$ and $c_2$, so it is conventional to
provide only the triple $(n,r,c_2)$.

\begin{theorem}
Suppose $X$ is an antipodal distance regular graph of diameter three with parameters
$(n,r,c_2)$ and let $A_i$ be the $i$-th distance matrix of $X$, for $i=1,2,3$.
Then the matrix
\[
W=I+xA_1+yA_2+zA_3
\]
is type-II if and only if 
\begin{enumerate}[(a)]
\item
$x=y$ and $W$ is a type-II matrix in the Bose-Mesner algebra of $rK_n$.
\item
$y=-x\inv$ and 
$x$ is a solution of a quadratic equation.
\item
$y\ne -x\inv$ and the possible values of $(x,y)$ are the points of intersection of 
two quartics in $x$ and $y$.
\end{enumerate}
\end{theorem}

\proof
We use $\theta$ and $\tau$ to denote eigenvalues of $X$ not equal to 
$-1$ or $n-1$.  Now $W$ is a type-II matrix if and only if 
the following system of equations are satisfied:
\begin{align}
(1-x-(r-1)y+(r-1)z) 
\left(1-\frac{1}{x}-\frac{(r-1)}{y}+\frac{(r-1)}{z}\right) &= nr, \label{Eqn1}\\
(1+\th x-\th y-z)(1+\th\ix-\th\iy-\iz) &=nr, \label{Eqn2}\\
(1+\tau x-\tau y-z) (1+\tau\ix-\tau \iy-\iz) &=nr. \label{Eqn3}
\end{align}
Subtracting \eqref{Eqn3} from \eqref{Eqn2} gives
\begin{equation}
\label{Eqn4}
(x-y)\iz+(\ix-\iy)z =(x-y)+(\ix-\iy)+(\theta+\tau)(x-y)(\ix-\iy).
\end{equation}
Adding $\th$ times this to \eqref{Eqn2} yields
\begin{equation}
\label{Eqn5}
\iz+z=-\theta \tau(x-y)(\ix-\iy)+2-nr.
\end{equation}

Solving \eqref{Eqn4} and \eqref{Eqn5} as two linear equations in $z$ and
$\iz$, we get
\begin{equation}
\label{Eqn6}
(x-y)\left((1+xy)z-\theta\tau(x-y)^2-(\theta+\tau)(x-y)+(nr-1)xy-1\right)=0.
\end{equation}
There are three cases.  First if $x=y$ we are lead to type-II matrices contained
in the Bose-Mesner algebra of $rK_n$ (including the Potts models).
Second, if $xy=-1$ then \eqref{Eqn6} yields a quadratic in $X:=x+x\inv$:
\begin{equation}
\label{Eqn7}
-\th\tau X^2-(\th+\tau)X-nr =0
\end{equation}
and \eqref{Eqn5} gives 
\begin{eqnarray}
z\inv + z 
&=& -\theta\tau X^2-nr+2\nonumber\\
&=& (\theta+\tau) X+2.\label{Eqn8}
\end{eqnarray}
Solving (\ref{Eqn1}) and (\ref{Eqn8}) as two linear equations in $z$ and
$z\inv$ gives 
\begin{equation*}
z = \frac{p(x)}{rx(x+1)(x-1)(r-1)}
\end{equation*}
and
\begin{equation*}
z\inv = \frac{-x^4p(x\inv)}{rx(x+1)(x-1)(r-1)}
\end{equation*}
where
\begin{eqnarray*}
p(x) &=& 
(r-1)(\theta+\tau+1)x^4 + (\theta+\tau+r-r\theta-r\tau)x^3+\\
&&
(3r\theta-r^2\tau-r^2\theta+3r-r\theta\tau+3r\tau-2-2\theta-2\tau-2r^2)x^2+\\
&&
(-r\theta-r\tau+3r+\theta+\tau-2r^2) x 
-(r-1)(r\theta+r\tau-1-\tau-\theta).
\end{eqnarray*}
Now these expressions for $z$ and $z\inv$ hold if and only if
\begin{equation*}
-x^4p(x\inv)p(x) = [rx(x+1)(x-1)(r-1)]^2
\end{equation*}
which gives a quartic in $X$.
Applying (\ref{Eqn7}) to this quartic, we can express
$X = x+x\inv$ in $r$, $\theta$, and $\tau$.
Hence $x$ is a solution of a quadratic equation.

Finally if $x\ne y$ or $-y\inv$, Equations~(\ref{Eqn4}) and (\ref{Eqn5}) 
are equivalent to
\begin{equation*}
z = \frac{1}{(1+xy)}
\left(\theta\tau(x-y)^2+(\theta+\tau)(x-y)-(nr-1)xy+1 \right),
\end{equation*}
and
\begin{equation*}
\iz = \frac{1}{xy(1+xy)}
\left(\theta\tau(x-y)^2-(\theta+\tau)(x-y)xy-(nr-1)xy+x^2y^2 \right).
\end{equation*}

Now substituting these two expressions into \eqref{Eqn1} gives a quartic in variables 
$x$ and $y$ while $z\iz=1$ gives another one.\qed

Note that $rK_n$ is a strongly regular graph, so the possible type-II matrices
are determined by the results of the previous section.

Calculations performed in Maple showed that
the resultant with respect to $x$ of the two quartics in case (c) is a non-zero
polynomial in $y$ of degree at most 30.  
By the elimination property of resultants \cite{MR2122859}, 
the resultant vanishes
at any common solution of the two quartics.  Hence these two quartics vanish
at no more than thirty values for $y$.
Similarly, the resultant with respect to $y$ of these two quartics is a non-zero
polynomial in $x$ of degree at most 30 and they vanish at no
more than thirty values for $x$.  Consequently there are finitely many
type-II matrices, up to scalar multiplication, in the Bose-Mesner algebra of an antipodal distance regular graph of diameter three.

As a final remark, it could be true that each Bose-Mesner algebra is equal to the set of all polynomials in some type-II matrix.  The results of the last two sections imply this is true
for schemes with at most two classes, and for antipodal schemes with three classes.
(Since we do not have strong evidence either way, we will not make any conjecture.)

\bibliography{adatyp2}
\bibliographystyle{acm}

\end{document}

%% file: Lmacs.tex
\newcommand\al{\alpha}
\newcommand\be{\beta}
\newcommand\de{\delta}

\newcommand\ga{\gamma}

\newcommand\la{\lambda}

\newcommand\om{\omega}
\newcommand\Om{\Omega}

\renewcommand\th{\theta} %--- Latex uses \th for a Norse character
\newcommand\Th{\Theta}

\newcommand\cA{{\mathcal A}}

\newcommand\cN{{\mathcal N}}

\newcommand\cx{{\mathbb C}}% complexes

\newcommand\re{{\mathbb R}}%reals

\newcommand\opk[1]{\mathop{\mathrm{#1}}\nolimits}

\newcommand\diff{\mathbin{\mkern-1.5mu\setminus\mkern-1.5mu}}% for \setminus

\newcommand\inv{^{-1}}
\newcommand\seq[3]{#1_{#2},\ldots,#1_{#3}}

\newcommand\pmat[1]{\begin{pmatrix} #1 \end{pmatrix}}

\newcommand\ip[2]{\langle#1,#2\rangle}
\newcommand\one{{\bf1}}
\newcommand\rk{\opk{rk}}
\newcommand\tr{\opk{tr}}